\documentclass[11pt,twoside]{article}
\usepackage{caption}

\usepackage[a6paper]{geometry}
\usepackage{amsmath}
\usepackage{resizegather}
\usepackage{epsfig,amsfonts,color}
\usepackage{amsmath}
\usepackage{graphicx}

\usepackage{amssymb, palatino, geometry,url}
\usepackage{algorithmic}
\usepackage[noresetcount,lined,boxed]{algorithm2e} 
\usepackage[colorlinks=true,linkcolor=blue,citecolor=blue,urlcolor=blue]{hyperref}
\usepackage{multirow}
\usepackage{fancyhdr}
\newcommand{\be}{\begin{equation}}
\newcommand{\ee}{\end{equation}}
\newcommand{\bea}{\begin{eqnarray}}
\newcommand{\eea}{\end{eqnarray}}

\newcommand{\bvec}{\left(\begin{array}{c}}
\newcommand{\evec}{\end{array}\right)}
\newcommand{\bsub}{\begin{subequations}}
\newcommand{\esub}{\end{subequations}}

\usepackage{amsthm} 

\usepackage{lineno}
\usepackage[sort]{cite} 
\usepackage{float} 
\usepackage{authblk}
\allowdisplaybreaks 
\usepackage{csvsimple}

\title{Mitigating Investment Risk Using Modular Technologies}
\author{Yue Shao, Yicheng Hu and Victor M. Zavala}
\date{}
\affil{Department of Chemical and Biological Engineering \\ University of Wisconsin-Madison, 1415 Engineering Dr., Madison, WI 53706, USA}

\begin{document}
\maketitle

\begin{abstract}
We study logistical investment flexibility provided by modular processing technologies for mitigating risk. Specifically, we propose a multi-stage stochastic programming formulation that determines optimal capacity expansion plans that mitigate demand uncertainty. The formulation accounts for multi-product dependencies between small/large units and for trade-offs between expected profit and risk. The formulation uses a cumulative risk measure to avoid time-consistency issues of traditional, per-stage risk-minimization formulations and we argue that this approach is more compatible with typical investment metrics such as the net present value. Case studies of different complexity are presented to illustrate the developments.  Our studies reveal that the Pareto frontier of a flexible setting (allowing for deployment of small units) dominates the Pareto frontier of an inflexible setting (allowing only for deployment of large units). Notably, this dominance is prevalent despite benefits arising from economies of scale of large processing units. 
\end{abstract}


\section{Introduction}
Modularization is a manufacturing trend that is being adopted in different industrial sectors such as power generation, data centers, and chemical processes \cite{Frivaldsky_2018, Berthelemy_2015, Dong_2009, Chakraborty_2009, DiPippo_1999}. Modularization enables technology size reduction and provides logistical flexibility to adapt to fast-changing markets and other externalities (e.g., climate, resource availability, policy)  \cite{Jaikumar_1986,Rajagopalan_1993}. For instance, decentralized power generation and storage systems are becoming increasingly attractive as climate change and adoption of renewable power disrupts markets and space-time demand patterns \cite{Heuberger_2017, Liu_2018,Shao_2019}.  Modular technologies can also be easily transported to different geographical locations to exploit changing market patterns and to enable the recovery of resources that are highly distributed and potentially short-lived \cite{allman2020dynamic,chen2019economies,Davie_2016}. We can interpret this ability as a form of {\em spatial-shifting} flexibility. This decentralized approach contrasts with the more traditional monolithic approach in which a large processing system is installed at a fixed location over its entire lifetime \cite{Zhao_2018}. This centralized approach involves investments that can reach billions of US dollars and face significant risk due to changing markets and climate, shortages of resources at a specific location (e.g., water), and changes in the policy landscape (e.g., carbon emissions).  As such, large central systems can face significant economic fallouts that investors might not be willing to tolerate. For instance, large ammonia production systems in the US have shut down due to low-cost supply from China, and large coal power plants are shutting down due to decreasing costs of renewable power and policy changes. Moreover, the mass deployment of small modular units facilitates experimentation, learning, and sharing of best practices that can reduce operational costs (compared to large facilities, in which experimentation is more difficult). On the downside, the flexibility provided by small modular systems often comes at the expense of increased investment and operational costs \cite{Rajagopalan_1993}.  Specifically, economies of scale benefit large systems due to the favorable scaling of throughput with equipment size \cite{Peters_1968}.  Due to complex trade-offs between costs and flexibility, industrial systems will likely evolve into a mixed state in which certain processing tasks are performed in small modular systems while others are performed in large centralized systems.  Identifying optimal investment strategies in such settings is complicated due to complex product interdependencies and uncertainties. 

A key observation driving this work is that modular systems provide logistical flexibility in investment size and timing that can be strategically exploited to mitigate risk. Specifically, expansion of production capacity in modular systems can proceed sequentially, which provides a mechanism to hedge against risk (we can interpret this as {\em temporal-shifting} flexibility). To give an example, the deployment of new power generators and transmission lines is subject to significant short-term and long-term uncertainties. Specifically, short-term fluctuations in demand and wind/solar supply can affect an optimal generation mix, and changes in fuel prices and policy can render entire technologies uneconomical \cite{Liu_2018}. Therefore, the progressive expansion of capacity using both large and small processing systems can help make and correct decisions and to better balance cost and risk.

In this work, we investigate investment flexibility provided by modular technologies; to do so, we propose a multi-product capacity expansion (CE) problem that exploits the availability of  technologies of different types and sizes to mitigate risk. Variants of the CE problem have been studied in different applications such as power generation, semiconductor manufacturing, railroad networks, and waste-to-energy systems \cite{Cardin_2015, Sun_2015, Takayuki_2018, Geng_2009}.  A cost-minimization CE problem that considers a single-product deterministic setting with installation decisions of a fixed-capacity facility was formulated in \cite{Luss_1979}. This formulation was extended to incorporate facilities with multiple capacities in \cite{ Luss_1983, Luss_1986}. Uncertainty in demand for a single-product cost-minimization CE problem was addressed by using a stochastic programming (SP) model in \cite{Murphy_1982, DapkLus_1984, Shiina_2003}. A stochastic CE formulation for planning investments in electricity generation, storage, and transmission investments over a long planning horizon was proposed in \cite{Liu_2018}.  These CE problem formulations use expected cost as an investment metric and thus do not control investment risk. Recently, a CE problem formulation that trades-off expected cost and risk was proposed in\cite{Zhao_2019}. Here, the conditional value-at-risk (CVaR) was used as a risk metric that is minimized at each stage. All the aforementioned formulations consider facilities that produce a single product; in a chemical process, however,  multi-product dependencies need to be captured. Specifically, a chemical manufacturing facility might involve processes that produce intermediate or final products and demands for such products might face different levels of uncertainty. Making investment decisions in a multi-product setting is a non-trivial problem. Capturing risk in time-dependent decision-making settings (such as CE) is also an active topic of research. For instance, time-consistency of per-stage risk minimization is an issue of concern. In the context of CE, time consistency indicates that, if an alternative A is riskier than alternative B at some time, then A should also be considered riskier than B at every prior time \cite{Boda_2006}.  Unfortunately, deriving SP formulations that achieve time-consistency is not straightforward.  Moreover, per-stage risk minimization is not necessarily a decision-making strategy that investors might follow; specifically, investors are typically  concerned with assessing risk of cumulative metrics such as the net present value (NPV). 

In this work, we propose a multi-product CE formulation to investigate flexibility brought by modularization for mitigating investment risk. Our framework is a multi-stage and multi-objective SP problem that captures demand product uncertainty and trade-offs between expected value and risk of the NPV. We provide case studies of different complexity to illustrate the developments. Our analysis reveals that the Pareto frontier of a flexible setting (allowing for deployment of units of various sizes) dominates the Pareto frontier of an inflexible setting (allowing only for deployment of large units). Our formulation also avoids difficulties associated with time-consistency issues of stage-wise risk-minimization formulations and we argue that is more compatible with more traditional investment strategies. 

\section{Problem Formulations}

In this section, we present CE formulations of different complexity (single-product/multi-product and deterministic/stochastic) in order to highlight different aspects of the problem. We begin our discussion by posing a couple of illustrative examples; this will help us introduce some key concepts that are essential in developing more complex CE formulations. 

\subsection{Problem Setting}

Consider the following deterministic CE setting: a decision-maker (investor) wants to progressively add capacity to a production system by installing technologies of different sizes (capacities). The resulting assembled system seeks to generate sufficient product to satisfy a time-dependent demand over a given planning horizon.  At each planning stage, the investor decides how many technologies (and associated capacities) it should install; if a technology is added at one stage, this will generate a product to satisfy the demand at the next stage (there is a deployment delay of one stage). Demand satisfaction generates revenue. We assume that an installed technology has to operate at full capacity; if the system production exceeds demand at a given time, the investor can decide to either store the excess product at a cost (and carry the product over to the next stage) or it can dispose of excess product at a cost. At the final stage, the system disposes of leftover excess product. The goal is to make an optimal CE plan over the horizon that maximizes NPV (accumulated cash flows over the horizon); in doing so, the investor is constrained by the capacities of the technologies available.  For simplicity, in this example, we assume that NPV is simply determined by the excess product (waste) at the end of the planning horizon and that there is no interest rate. 

\begin{figure}[hbt!]
\centering
    \includegraphics[scale=0.4]{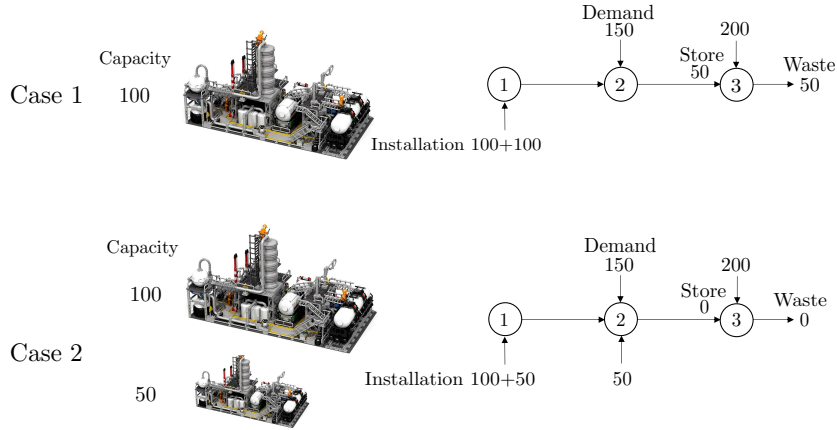}
    \caption{Illustrative example of the single-product deterministic capacity expansion setting.}
    \label{fig:illus_det}
\end{figure} 

\begin{figure}[hbt!]
\centering
    \includegraphics[scale=0.4]{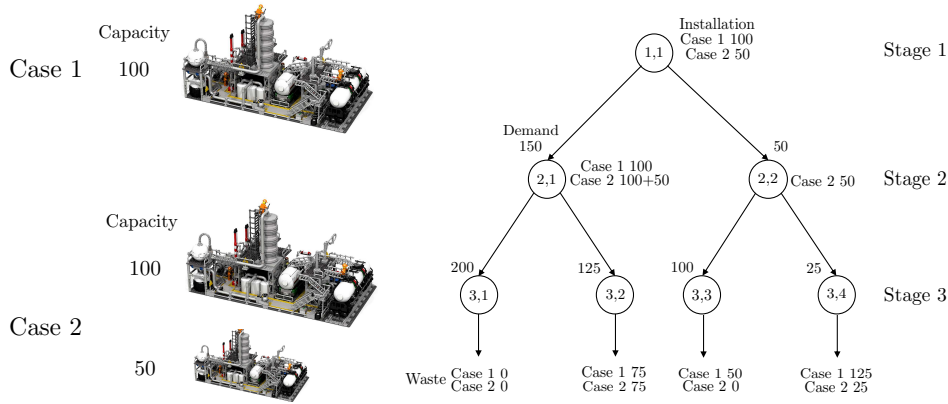}
    \caption{Illustrative example of the single-product stochastic capacity expansion setting.}
    \label{fig:illus_sto}
\end{figure}

We illustrate this decision-making setting in Figure \ref{fig:illus_det}; here, we would like to make decisions on how much capacity to install at Stage 1 and Stage 2 to minimize waste at Stage 3. In Case 1, only large technologies are available (with a capacity of 100 units); to satisfy future demands, it is decided to install 2 units of this large technology at Stage 1. Since the demand at Stage 2 is 150, it is required to shift excess production to Stage 3. Moreover, since the demand at Stage 3 is 200, it is necessary to dispose of 50 units of excess product.  In Case 2, a large technology and a small  technology are available; this opens the possibility of an investment strategy in which we install a large unit  (size 100) and a small unit (size 50) at Stage 1 and add a small unit at Stage 2.  This strategy prevents wasting material at Stage 3 and highlights the flexibility provided by the availability of small technologies. Note that, in this setting, the demands are time-dependent but are assumed to be known at the moment of decision (deterministic setting). 

The CE setting can be extended to account for uncertainty in the demands (stochastic setting). Here, demand uncertainty is represented in the form of possible scenarios. An illustrative example of this setting is shown in Figure \ref{fig:illus_sto}. We would like to make installation decisions at each stage and scenario (here, we consider two possible scenarios per stage). Stages and scenarios are represented as a decision tree and each node is associated with a different demand scenario. Installation decisions are shown next to the node and wasted amounts are shown exiting the nodes at the last stage. In Case 1 (only large technologies are available), we decide to install a large technology in Stage 1; in Stage 2, we can decide to install a large technology in scenario 1 (high demand) or no technology in Scenario 2 (low demand). This investment strategy results in four scenarios of waste product in Stage 3 (10,75,50,125). Assuming that these scenarios have equal probability (1/4), the expected value of the waste is 65 and the standard deviation (typical measure of risk) is 48.  In Case 2 (large and small technologies available), we install a small unit in Stage 1; in Scenario 1 in Stage 2 we install a large technology (to satisfy the large demand) and in Scenario 2 we install another small technology (to satisfy the small demand). This investment strategy results in four scenarios of waste excess product in Stage 3 (0,75,0,25). This gives a mean waste of 25 and a risk of 35. We can thus see that adding the possibility of installing small units reduces expected waste and risk.  
\\

Risk can be measured in different ways; in the previous setting, we computed the risk at Stage 3 (last stage) but we could have also computed the risk at Stage 2 and we could have added this to the risk of Stage 3 (add risks for all stages) to determine the best strategy. This highlights issues that one may encounter when measuring risk in a multi-stage decision-making setting. Specifically, risk can vary over  time and one might or might not be interested in shaping risk over time. This is similar in spirit to how investors think about cash flows; typically, investors are not necessarily interested in the temporal behavior of cash flows but want to aggregate cash flows in a single metric (e.g., NPV). Following this reasoning, in this work, we will compute NPV for every branch in the tree and compute the associated risk.
\\

The CE problem can be further extended to a multi-product setting in which a system can produce multiple intermediate or final products. Intermediate products generate interdependencies between possible technologies (i.e., technology can take intermediate products obtained from another technology as raw materials).  Multi-product dependencies make the problem significantly more complicated and we will see that,  in such a setting, investment flexibility provided by small units becomes particularly relevant. We now proceed to formulate single-product deterministic and stochastic CE problems and we then proceed to extend this to a multi-product setting. 

\begin{figure}[hbt!]
\centering
    \includegraphics[scale=0.6]{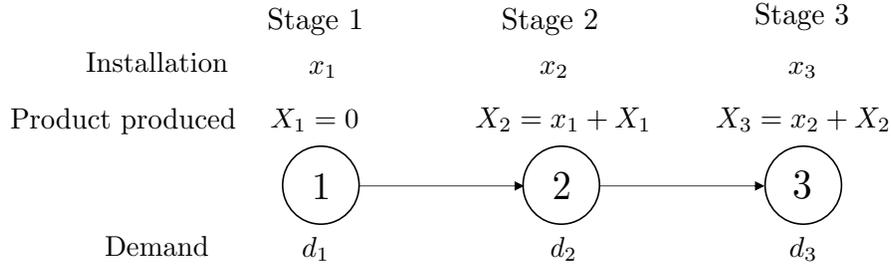}
    \caption{Tree representation of planning horizon in deterministic case}
    \label{fig:deter}
\end{figure} 

\subsection{Single-Product, Deterministic Setting}

Consider the decision-making setting shown in Figure \ref{fig:deter}. We consider a planning horizon comprising a set of stages  $\mathcal{T}=\{1,2,\dots, T\}$ with cardinality $|\mathcal{T}|:=T$.  The time-dependent product demand is given by  $d_t, t \in \mathcal{T}$. Investment decisions are made at stages $t \in \{1,2,\dots,T-1\}$  and we thus define the decision stages $\mathcal{D}=\{1,2,\dots,T-1\}$ with cardinality $|\mathcal{D}|:=T-1$. In a deterministic setting, the  planning horizon is a linear graph (a tree) in which each node represents a stage.  As such, for each node $t$, we define a parent node $a_t\in \mathcal{T}$  (in this case we have $a_t=t-1$). The root node $t=1$ does not have a parent node and thus $a_1 = \emptyset$. 
\\

The investor has a list of of possible technology choices that can be installed at each stage. Each choice has a different capacity and associated installation cost (which capture economies of scale). We define the set of capacities as $\mathcal{B}=\{B_1, B_2, \dots, B_N\} \in \mathbb{Z}^N_+$ and the set of associated costs as $\mathcal{C}=\{C_1,C_2,\dots,C_N\} \in \mathbb{R}^N_+$, both with same cardinality $|\mathcal{B}|=|\mathcal{C}|:=N$. For convenience, we also define a set of choice indexes $\mathcal{F} = \{1, 2, \dots, N\}$.  To capture economies of scale, it is typical to assume that costs follow a {\em 2/3 rule} and thus: 
\begin{align}
\left(\frac{B_{k}}{B_{k'}}\right) = \left(\frac{C_{k}}{C_{k'}}\right)^{\frac{3}{2}}, \hspace{0.2cm} k, k' \in \mathcal{F}
\end{align}
where $B_{k}, B_{k'} \in \mathcal{B}$ are the ${k}^{th}$ and ${k'}^{th}$ capacity choices and $C_{k}, C_{k'} \in \mathcal{C}$ are the installation costs. 
\\

Product storage comes at a cost $\rho_s \in \mathbb{R}+$ and we define a maximum storage capacity $\bar{s} \in \mathbb{Z}_+$. Disposal of excess product comes at a cost $\rho_w$. We define a variable $s_{t} \in \mathbb{Z}_+, t \in \mathcal{T}$ to capture the amount of storage at stage $t$. We set $s_1=0$ and $s_n=0$ (any excess product is regarded as waste at the final stage). We define the integer variable $w_{t} \in \mathbb{Z}_+, t \in \mathcal{T}$ to represent the waste generated at each stage. We assume $w_1=0$ (waste is generated at the end of each stage). The investor has a choice to deal with any excess product; either to dispose of the product or to store it (shift it to the next stage). To capture installation {\em delays}, we assume that capacity installed at stage $t$ generate  production, storage, disposal and sales of products at stage $t+1$.
\\

We define integer variables $u_{t,k} \in \mathbb Z_+, t \in \mathcal{D}, k \in \mathcal{F}$;  here, $u_{t,k}$ is the number of technologies of type $k\in \mathcal{B}$ installed at stage $t\in \mathcal{D}$. The total capacity installed at time $t$ is thus:
\begin{equation}\label{eq:d_x}
x_{t}=\sum_{k \in \mathcal{F}}u_{t,k}B_k, \hspace{0.2cm} t \in \mathcal{D}
\end{equation} 
and the total installation cost at time $t$ is:
\begin{equation}\label{eq:d_y}
y_{t} =\sum_{k \in \mathcal{F}}u_{t,k}C_k, \hspace{0.2cm} t \in \mathcal{D}
\end{equation}
We define variable $X_t, t \in \mathcal{T}$ to represent the total amount of product generated at stage $t$; $X_t$ is the cumulative installed capacity up to stage $t$ and follows the dynamic evolution: 
\begin{align}
\label{eq:d_cu_x}
&\; X_t=X_{a_t}+x_{a_t}, \hspace{0.2cm} t \in \mathcal{T}
\end{align}
We set the initial production as $X_1=0$ and recall that $a_t=t-1$. In the proposed setting, we can install more than one technology at each time but we limit the total final installed capacity $X_T$ by using the upper bound $\bar{x}\in \mathbb{Z}_+$. 
\\

We use the following constraint to ensure that demand is satisfies at each stage:
\begin{align}
& \; 0 \leq s_{a_t}+X_t-s_{t}-w_t \leq d_{t}, \hspace{0.2cm} t \in \mathcal{T}\label{eq:d_sin_ba}
\end{align}
We define a production cost as $\rho_p\in \mathbb{R}_+$ and selling price as $\pi_p \in \mathbb{R}_+$. Under these definitions, the cost at stage $t$ (denoted as $q_{t}\in \mathbb{R}$) can be expressed as:
\begin{subequations}\label{eq:d_q}
\begin{align}
&\; q_{t}=y_{t}+\rho_pX_t+\rho_{s}s_{t}+\rho_{w}w_{t}, \hspace{0.2cm} t \in \mathcal{T}.
\end{align}
\end{subequations}
and the revenue at stage $t$ (denoted as $r_{t}\in \mathbb{R}$) is expressed as: 
\begin{align}\label{eq:d_v}
&\; r_{t}=\pi_p(X_t+s_{a_t}-s_t-w_t), \hspace{0.2cm} t \in \mathcal{T}.
\end{align}
Note that $s_0,X_1,s_1,w_1=0$ and thus $r_1=0$.  
\\

We consider a CE formulation that maximizes the NPV of the project; to do so, we define an interest rate $\gamma \in [0,1]$ that is used to discount any future cash flow and we define the discount factor $\beta_t=1/(1+\gamma)^{t-1}$.  We define the discounted profit (cash flow) achieved at stage $t$ as $v_{t}$ and the cumulative profit upto stage $t$ as $V_{t}$. These quantities are computed as:
\begin{subequations}
\begin{align}
&\; v_{t}=\beta_t\cdot (r_{t}-q_{t}), \hspace{0.2cm} t \in \mathcal{T}\\
&\; V_{t}=V_{a_{t}}+v_{t}, \hspace{0.2cm} t \in \mathcal{T}.
\end{align}
\end{subequations}
With this, the NPV is given by $V_T=\sum_{t \in \mathcal{T}}  v_t$. 
\\

In summary, the CE problem is a mixed-integer linear program (MILP) of the form:
\begin{subequations}\label{eq:d_single_pro}
\begin{align}
\max_{u,s}\; &\; V_T\\
\textrm{s.t.} 
& \;x_{t}=\sum_{k \in \mathcal{F}}u_{t,k}B_k, \hspace{0.2cm} t \in \mathcal{D}\\
& \; y_{t} =\sum_{k \in \mathcal{F}}u_{t,k}C_k, \hspace{0.2cm} t \in \mathcal{D}\\
& \; X_t=X_{a_t}+x_{a_t}, \hspace{0.2cm} t \in \mathcal{T}\\
& \; q_{t}=y_{t}+\rho_pX_t+\rho_{s}s_{t}+\rho_{w}w_{t}, \hspace{0.2cm} t \in \mathcal{T}\\
& \; r_{t}=\pi_p(X_t+s_{a_t}-s_t-w_t), \hspace{0.2cm} t \in \mathcal{T}\\
&\; v_{t}=\beta_t\cdot (r_{t}-q_{t}), \hspace{0.2cm} t \in \mathcal{T}\\
&\; V_{t}=V_{a_{t}}+v_{t}, \hspace{0.2cm} t \in \mathcal{T}\\
& \; 0 \leq s_{a_t}+X_t-s_{t}-w_t \leq d_{t}, \hspace{0.2cm} t \in \mathcal{T}\\
& \; 0 \leq s_{t} \leq \bar{s}, \hspace{0.2cm} t \in \mathcal{D}\\
& \; X_T \leq \bar{x}\\
& \; u_{t,k} \in \mathbb Z_+, \hspace{0.2cm} t \in \mathcal{D}, k \in \mathcal{F}\\
& \; s_{t} \in \mathbb Z_+, \hspace{0.2cm} t \in \mathcal{T}.
\end{align}
\end{subequations}
The NPV metric accumulates all the cash flows $v_t,\,t\in \mathcal{T}$ to the initial stage $t=1$ and this accounts for time value of money. If we set $\gamma=0$, we obtain $\beta_t=1$ and the CE problem maximizes the cumulative cash flows over the planning horizon (the total profit).  As we discuss next, the NPV is a convenient metric that allows us to summarize random cash flows that arise in settings that face uncertainty. 

 \begin{figure}[hbt!]
\centering
    \includegraphics[scale=0.6]{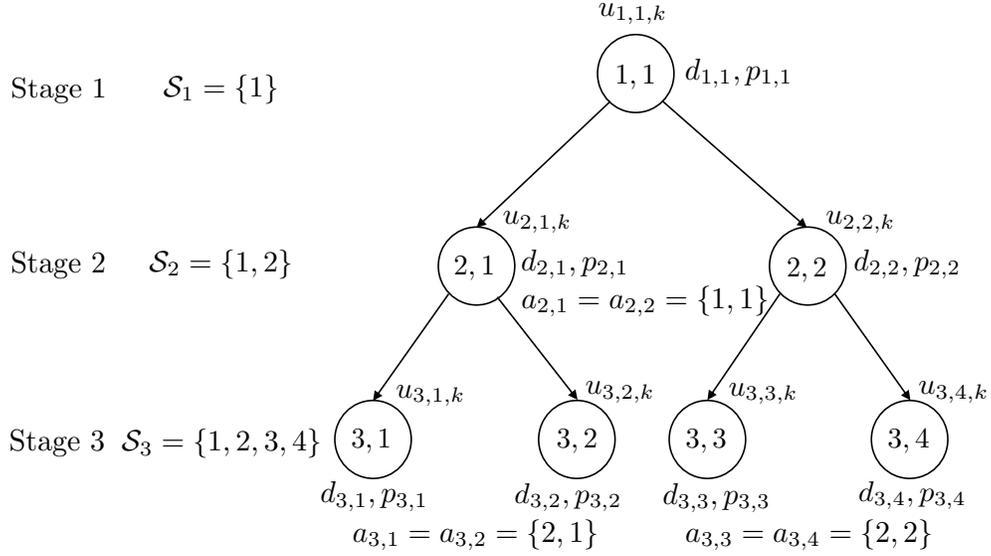}
    \caption{Tree representation of planning horizon in stochastic case}
    \label{fig:single_for}
\end{figure} 

\subsection{Single-Product, Stochastic Setting}

\begin{figure}[!ht]
\centering
    \includegraphics[scale=0.5]{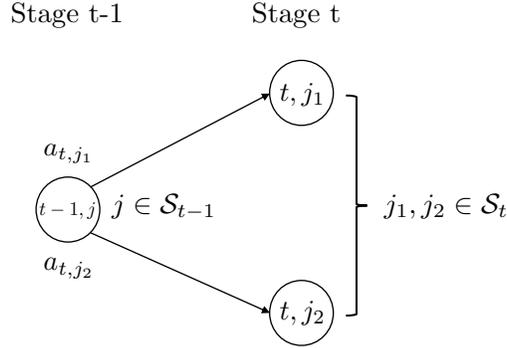}
    \caption{Schematic of parent-node notation. Here, node $\{t-1, j\}$ is the ancestor of $\{t-1, j_1'\}$ and $\{t-1, j\}$ and thus $a_{t,j_1}=a_{t,j_2}=\{t-1, j\}$.} 
    \label{fig:parent}
\end{figure}

We now extend the CE problem to a stochastic setting; this formulation allows us to explore trade-offs between expected profit and risk. The stochastic setting is illustrated in Figure \ref{fig:single_for}. We define the set of scenarios at each stage $t \in \mathcal{T}$ as $\mathcal{S}_t=\{1,2,\dots S_t\}$ with cardinality $|\mathcal{S}_t|:=S_t$. Each scenario is represented as a node in a tree; the number of levels in the tree is given by the number of stages. We define parent node $a_{t,j}, t \in \mathcal{T}, j \in \mathcal{S}_t$ as the parent stage and scenario that node $\{t,j\}$ emanates from. For example, if scenario $\{t,j'_1\}$ is generated from scenario $j$ at stage $t-1$, then the parent node is $a_{t,j'_1}=\{t-1,j\}$ (see Figure \ref{fig:parent}). The scenario set of the root node is a singleton $\mathcal{S}_1 = \{1\}$ and the parent of the root node is empty and thus $a_{1,1} = \emptyset$.

The demand is a discrete random variable; the realization of this variable at time $t$ and scenario $j$ is denoted as $d_{t,j}, t \in \mathcal{T}, j \in \mathcal{S}_t$.  The probability of realization $d_{t,j}$ is represented as $p_{t,j}\in [0,1]$.  For each stage $t\in \mathcal{T}$, these probabilities satisfy $\sum_{j \in \mathcal{S}_t} p_{t,j}=1$. It is important to highlight that these are {\em joint} probabilities that capture the {\em history of events} leading to node $\{t,j\}$. In other words, joint probabilities are node probabilites (conditional probabilities are arc/edge probabilities and marginal probabilities ignore history). For example, in Figure \ref{fig:single_for}, $p_{3,1}$ is the probability of node $\{3,1\}$ corresponding to the demand event sequence $d_{1,1}, d_{2,1}$ and $d_{3,1}$ and thus: 
\begin{subequations}\label{eq:prob}
\begin{align}
p_{3,1} & \; =\mathbb{P}(\{1,1\}, \{2,1\}, \{3,1\})\\
& \; =\mathbb{P}(\{3,1\}|\{1,1\}, \{2,1\})\cdot \mathbb{P}(\{1,1\}, \{2,1\})\\
& \; =\mathbb{P}(\{3,1\}|\{1,1\}, \{2,1\})\cdot p_{2,1}\\
& \; =\mathbb{P}(\{3,1\}|\{1,1\}, \{2,1\})\cdot \mathbb{P}(\{2,1\}|\{1,1\})\cdot p_{1,1}.
\end{align}
\end{subequations} 
where $\mathbb{P}(\{3,1\}|\{1,1\}, \{2,1\})$ is the conditional probability of event $\{3,1\}$ given that $\{1,1\}$ and $\{2,1\}$ have been realized. We thus have that probability $p_{t,j}$ carries information of its ancestor nodes (i.e., $p_{3,1}$ carries information of $p_{2,1}$ and $p_{1,1}$). 
\\

 We define an integer variable $u_{t,j,k} \in \mathbb Z_+, t \in \mathcal{D}, j \in \mathcal{S}_t, k \in \mathcal{F}$ to represent the number of technologies of type $k\in\mathcal{B}$ installed at stage $t\in \mathcal{T}$ and at scenario $j\in \mathcal{S}_t$. The total capacity installed at time $t$ and scenario $j$ is: 
\begin{equation}\label{eq:x}
x_{t,j}=\sum_{k \in \mathcal{F}}u_{t,j,k}B_k, \hspace{0.2cm} t \in \mathcal{D}, j \in \mathcal{S}_t,
\end{equation} 
and total installation cost at time $t$  and scenario $j$ is:
\begin{equation}\label{eq:y}
y_{t,j} =\sum_{k \in \mathcal{F}}u_{t,j,k}C_k, \hspace{0.2cm} t \in \mathcal{D}, \hspace{0.2cm} j \in \mathcal{S}_t.
\end{equation}
We use the integer variable $s_{t,j} \in \mathbb Z_+, t \in \mathcal{T}, j \in \mathcal{S}_t$ to denote the amount of storage at stage $t$ in scenario $j$. Similar to the deterministic case, the storage for the first and last stage are assumed to be zero. We define the waste variable as $w_{t,j} \in \mathbb Z_+, t \in \mathcal{T}, j \in \mathcal{S}_t$ and we assume that the waste for the first stage zero. 
\\

We define a variable $X_{t,j}, t \in \mathcal{T}, j \in \mathcal{S}_t$ as the total amount of product produced at stage $t$ and scenario $j$. Here, $X_{t,j}$ is interpreted as the cumulative installed capacity up to stage and scenario $a_{t,j}$, expressed in \eqref{eq:cu_x}. The total installed capacity at stage one is $X_{1,j}=0, j \in \mathcal{S}_1$.
\begin{align}
\label{eq:cu_x}
\; X_{t,j}=X_{a_{t,j}}+x_{a_{t,j}}, \hspace{0.2cm} t \in \mathcal{T}, j \in \mathcal{S}_t
\end{align}
With the above definitions, we can define the undiscounted cost at stage $t$ for scenario $j$ as $q_{t,j}$, and is given by:
\begin{subequations}\label{eq:q}
\begin{align}
&\; q_{t,j}=y_{t,j}+\rho_pX_{t,j}+\rho_{s}s_{t,j}+\rho_{w}w_{t,j}, \hspace{0.2cm} t \in \mathcal{T}\setminus\{T\}, j \in \mathcal{S}_t.
\end{align}
\end{subequations}
The undiscounted revenue at stage $t$ scenario $j$ is denoted as $r_{t,j}$ and can be expressed as
\begin{subequations}\label{eq:v}
\begin{align}
&\; r_{t,j}=\pi_p(X_{t,j}+s_{a_{t,j}}-s_{t,j}-w_{t,j}), \hspace{0.2cm} t \in \mathcal{T}, j \in \mathcal{S}_t.
\end{align}
\end{subequations}

Our goal is to maximize the expected NPV and to minimize its risk.  To model these quantities, we introduce variable $v_{t,j}, t \in \mathcal{T}, j \in \mathcal{S}_t$ that denotes the cash flow (profit) achieved in stage $t$ and scenario $j$. We also define the cumulative variable $V_{t,j}, t \in \mathcal{T}, j \in \mathcal{S}_t$ to denote the cumulative profit up to stage $t$ and scenario $j$. Using the notation proposed, these quantities can be computed using a form that is analogous to the deterministic case:
\begin{subequations}
\begin{align}
&\; v_{t,j}=\beta_t \cdot (r_{t,j}-q_{t,j}), \hspace{0.2cm} t \in \mathcal{T}, j \in \mathcal{S}_t\\
&\; V_{t,j}=V_{a_{t,j}}+v_{t,j}, \hspace{0.2cm} t \in \mathcal{T}, j \in \mathcal{S}_t.
\end{align}
\end{subequations}
The NPV is the total accumulated cash flow and is given by $V_{T,j},\, j\in \mathcal{S}_T$. We note that this is a random quantity and that each realization correspond to a branch of the scenario tree connecting the root node $\{1,1\}$ to the final nodes $\{T,j\}$ with $j\in \mathcal{S}_T$.  The  NPV thus summarizes information of the entire project and captures probabilities of the different paths that the project can take. When the interest rate is zero, the NPV of a given path is the total profit of the project for such path. 
\\

The expected NPV is given by:
\begin{align}
&\; \mathcal{E}=\sum_{j \in \mathcal{S}_T}p_{T,j}V_{T,j},
\end{align}
and its risk is measured by using the mean deviation:
\begin{align}
&\; \mathcal{R}=\sum_{j \in \mathcal{S}_T}p_{T,j}|V_{T,j}-\mathcal{E}|.
\end{align}
Alternative risk metrics can be used; here, we provide the mean deviation as this is a coherent risk measure that is easy to interpret. 
\\

The CE problem can be cast as the following stochastic, multistage, multiobjective optimization (SMMO) problem:
\begin{subequations}\label{eq:single_risk}
\begin{align}
\max_{u,s}\; &\; \{\mathcal{E},-\mathcal{R}\}\\
\textrm{s.t.}& \; x_{t,j}=\sum_{k \in \mathcal{F}}u_{t,j,k}B_k, \hspace{0.2cm} t \in \mathcal{D}, j \in \mathcal{S}_t\\
& \; y_{t,j} =\sum_{k \in \mathcal{F}}u_{t,j,k}C_k, \hspace{0.2cm} t \in \mathcal{D}, \hspace{0.2cm} j \in \mathcal{S}_t\\
& \; X_{t,j}=X_{a_{t,j}}+x_{a_{t,j}}, \hspace{0.2cm} t \in \mathcal{T}, j \in \mathcal{S}_t\\
&\; q_{t,j}=y_{t,j}+\rho_pX_{t,j}+\rho_{s}s_{t,j}+\rho_{w}w_{t,j}, \hspace{0.2cm} t \in \mathcal{T}, j \in \mathcal{S}_t\\
& \; r_{t,j}=\pi_p(X_{t,j}+s_{a_{t,j}}-s_{t,j}-w_{t,j}), \hspace{0.2cm} t \in \mathcal{T}, j \in \mathcal{S}_t\\
&\; v_{t,j}=\beta_t(r_{t,j}-q_{t,j}),  \hspace{0.2cm} t \in \mathcal{T}, j \in \mathcal{S}_t\\
&\; V_{t,j}=V_{a_{t,j}}+v_{t,j}, \hspace{0.2cm} t \in \mathcal{T}, j \in \mathcal{S}_t\\
&\; \mathcal{E}=\sum_{j \in \mathcal{S}_T}p_{T,j}V_{T,j}\\
&\; \mathcal{R}=\sum_{j \in \mathcal{S}_T}p_{T,j}|V_{T,j}-\mathcal{E}|\\
& \; 0 \leq s_{a_{t,j}}+X_{t,j}-s_{t,j}-w_{t,j} \leq d_{t,j}, \hspace{0.2cm} t \in \mathcal{T} j \in \mathcal{S}_t\\\label{eq:risk_con}
& \; 0 \leq s_{t,j} \leq \bar{s}, \hspace{0.2cm} t \in \mathcal{D}, j \in \mathcal{S}_t\\
& \; X_{T,j} \leq \bar{x}, \hspace{0.2cm} j \in \mathcal{S}_T\\
& \; u_{t,j,k} \in \mathbb Z_+, \hspace{0.2cm} t \in \mathcal{D}, j \in \mathcal{S}_t, k \in \mathcal{F}\\
& \; s_{t,j} \in \mathbb Z_+, \hspace{0.2cm} t \in \mathcal{T}, j \in \mathcal{S}_t
\end{align}
\end{subequations}
The Pareto solutions of this problem are found by using an $\epsilon$-constrained method. It is important to highlight that the SMMO problem does not seek to optimize the conditional expectation and risk at every time (as in traditional multi-stage SP formulations). Instead, the SMOO problem optimizes the {\em joint} expectation and risk (over the entire planning time). This formulation thus avoids ambiguity issues associated with time consistency of conditional risk evaluation encountered in traditional formulations. Another way to think about this difference is that our formulation first determines the accumulated cash flow over all stages and then optimizes its risk, while a traditional formulation determines the risk of the cash flow at each stage and then optimizes the accumulated risk over all stages.  

\subsection{Multi-Product, Stochastic Setting}
We can conveniently extend the previous formulation to a multi-product setting. Assume that the investor now has a choice of producing multiple products $\mathcal{I}=\{i_1, i_2, \dots, i_I\}$. We use $\alpha_{i,i'}, i,i'\in \mathcal{I}$ to represent the inter-dependencies between product $i$ and $i'$. Specifically, $\alpha_{i,i'}$ denotes the units of product $i'$ required to produce $i$. Note that $\alpha_{i,i}=0, i \in \mathcal{I}$. For each product $i \in \mathcal{I}$, we define a set of technologies that can produce it; these technologies have capacities $\mathcal{B}^i$ and installation costs $\mathcal{C}^i$. Also, for each product $i \in \mathcal{I}$, we define a storage cost, waste disposal cost, operational cost, and selling price as $\rho_s^i$, $\rho_w^i$, $\rho_p^i$, and $\pi_p^i$. We also define a capacity limit for product $i$ as $\bar{x}^i,  i \in \mathcal{I}$, and storage limit for product $i$ as $\bar{s}^i,  i \in \mathcal{I}$. We define the demand for product $i$ at stage $t$ and scenario $j$ as $d_{t,j}^i, t \in \mathcal{T}, j \in \mathcal{S}_t, i \in \mathcal{I}$. Our decision variables are the number of technologies with capacity $B_k^i$ to be installed from the capacities list at stage $t$ and scenario $j$ for product $i$ and these are modeled using the integer variables  $u_{t,j,k}^i \in \mathbb Z_{+}, t \in \mathcal{D}, j \in \mathcal{S}_t, k \in \mathcal{F}^i, i \in \mathcal{I}$. 
\\

The total capacity installed at stage $t$ and scenario $j$ for product $i$ is denoted as $x_{t,j}^i, t \in \mathcal{D}, j \in \mathcal{S}_t, i \in \mathcal{I}$. The total installation cost at stage $t$ and scenario $j$ for product $i$ is denoted as $y_{t,j}^i, t \in \mathcal{D}, j \in \mathcal{S}_t, i \in \mathcal{I}$. These quantities are computed as:
\begin{subequations}
\begin{align}
&\; x_{t,j}^i=\sum_{k \in \mathcal{F}^i}u_{t,j,k}^iB_k^i, \hspace{0.2cm}  t \in \mathcal{D}, j \in \mathcal{S}_t, i \in \mathcal{I}\\
&\; y_{t,j}^i=\sum_{k \in \mathcal{F}^i}u_{t,j,k}^iC_k^i, \hspace{0.2cm}   t \in \mathcal{D}, j \in \mathcal{S}_t, i \in \mathcal{I}
\end{align}
\end{subequations} 
The amount of storage at stage $t$ and scenario $j$ for product $i$ is defined as $s_{t,j}^i \in \mathbb Z_+, t \in \mathcal{T}, j \in \mathcal{S}_t, i \in \mathcal{I}$, and the amount of product disposed at stage $t$ and scenario $j$ for product $i$ is $w_{t,j}^i \in \mathbb Z_+, t \in \mathcal{T}, j \in \mathcal{S}_t, i \in \mathcal{I}$. The total production at stage $t$ and scenario $j$ for product $i$ is  denoted as $X_{t,j}^i,t \in \mathcal{D}, j \in \mathcal{S}_t, i \in \mathcal{I}$. We also incorporate the cumulative installation cost occurred along the path to time $t$ and scenario $j$ for product $i$ and denote this as $Y_{t,j}^i, t \in \mathcal{D}, j \in \mathcal{S}_t, i \in \mathcal{I}$. These quantities are computed as:
\begin{subequations}
\label{eq:mul_cu}
\begin{align}
&\; X_{t,j}^i=X_{a_{t,j}}^i+x_{a_{t,j}}^i, \hspace{0.2cm}  t \in \mathcal{T}, j \in \mathcal{S}_t, i \in \mathcal{I}\\
&\; Y_{t,j}^i=Y_{a_{t,j}}^i+y_{t,j}^i, \hspace{0.2cm}  t \in \mathcal{D}, j \in \mathcal{S}_t, i \in \mathcal{I}
\end{align}
\end{subequations}
The cost incurred at stage $t$ and scenario $j$, is $q_{t,j}, j \in \mathcal{S}_T$ and is computed as:
\begin{align}
&\; q_{t,j}=\sum_{i \in \mathcal{I}}y_{t,j}^i+\rho_{p}^iX_{t,j}^i+\rho_{s}^is_{t,j}^i+\rho_{w}^iw_{t,j}^i, \hspace{0.2cm}  t \in \mathcal{T}, j \in \mathcal{S}_t.
\end{align}
The profit incurred at stage $t$ and scenario $j$ is $r_{t,j}, t \in \mathcal{D}, j \in \mathcal{S}_t$ and is computed as:
\begin{align}
&\; r_{t,j}=\sum_{i \in \mathcal{I}}\pi_p^i(X_{t,j}^i+s_{a_{t,j}}^i-s_{t,j}^i-w_{t,j}^i-\sum_{i' \in \mathcal{I}}X_{t,j}^{i'}\alpha_{i',i}), \hspace{0.2cm}  t \in \mathcal{T}, j \in \mathcal{S}_t.
\end{align}
Under these definitions, we can define the rest of the quantities for cash flow, cumulative cash flow, and NPV in the same way that we did for the single-product case. This gives the SMMO problem: 

\begin{subequations}
\label{eq:multi_risk}
\begin{align}
\max_{u,s}\; &\;\{\mathcal{E},-\mathcal{R}\}\\
\textrm{s.t.}&\; x_{t,j}^i=\sum_{k \in \mathcal{F}^i}u_{t,j,k}^iB_k^i, \hspace{0.2cm}  t \in \mathcal{D}, j \in \mathcal{S}_t, i \in \mathcal{I}\\
&\; y_{t,j}^i=\sum_{k \in \mathcal{F}^i}u_{t,j,k}^iC_k^i, \hspace{0.2cm}   t \in \mathcal{D}, j \in \mathcal{S}_t, i \in \mathcal{I}\\
&\; X_{t,j}^i=X_{a_{t,j}}^i+x_{a_{t,j}}^i, \hspace{0.2cm}  t \in \mathcal{T}, j \in \mathcal{S}_t, i \in \mathcal{I}\\
&\; Y_{t,j}^i=Y_{a_{t,j}}^i+y_{t,j}^i, \hspace{0.2cm}  t \in \mathcal{D}, j \in \mathcal{S}_t, i \in \mathcal{I}\\
&\; q_{t,j}=\sum_{i \in \mathcal{I}}y_{t,j}^i+\rho_{p}^iX_{t,j}^i+\rho_{s}^is_{t,j}^i+\rho_{w}^iw_{t,j}^i, \hspace{0.2cm}  t \in \mathcal{D}, j \in \mathcal{S}_t\\
&\; r_{t,j}=\sum_{i \in \mathcal{I}}\pi_p^i(X_{t,j}^i+s_{a_{t,j}}^i-s_{t,j}^i-w_{t,j}^i-\sum_{i' \in \mathcal{I}}X_{t,j}^{i'}\alpha_{i',i}), \hspace{0.2cm}  t \in \mathcal{T}, j \in \mathcal{S}_t\\
&\; v_{t,j}=\beta_t(r_{t,j}-q_{t,j}),  \hspace{0.2cm} t \in \mathcal{T}, j \in \mathcal{S}_t\\
&\; V_{t,j}=V_{a_{t,j}}+v_{t,j}, \hspace{0.2cm} t \in \mathcal{T}, j \in \mathcal{S}_t\\
&\; \mathcal{E}=\sum_{j \in \mathcal{S}_T}p_{T,j}V_{T,j}\\
&\; \mathcal{R}=\sum_{j \in \mathcal{S}_T}p_{T,j}|V_{T,j}-\mathcal{E}|\\
& \; 0 \leq s_{a_{t,j}}^i+X_{t,j}^i-s_{t,j}^i-w_{t,j}^i-\sum_{i \in \mathcal{I}}X_{t,j}^{i'}\alpha_{i',i} \leq d_{t,j}^i, \hspace{0.2cm}  t \in \mathcal{T}, j \in \mathcal{S}_t, i,i' \in \mathcal{I}\\
& \; \alpha_{i,i'}X_{t,j}^i \leq s_{a_{t,j}}^{i'}+X_{t,j}^{i'}, \hspace{0.2cm}  t \in \mathcal{D}, j \in \mathcal{S}_t, i,i' \in \mathcal{I}\\
& \; 0 \leq s_{t,j}^i \leq \bar{s}^i, \hspace{0.2cm}  t \in \mathcal{D}, j \in \mathcal{S}_t, i \in \mathcal{I}\\
& \; X_{T,j}^i \leq \bar{x}^i, \hspace{0.2cm}  j \in \mathcal{S}_T, i \in \mathcal{I}\\
& \; \sum_{i \in \mathcal{I}} Y_{T,j}^i \leq \bar{y}, \hspace{0.2cm}  j \in \mathcal{S}_T\\
& \; u_{t,j,k}^i \in \mathbb Z_+, \hspace{0.2cm}  i \in \mathcal{D}, j \in \mathcal{S}_t, k \in \mathcal{F}^i, i \in \mathcal{I}\\
& \; s_{t,j}^i \in \mathbb Z_+, \hspace{0.2cm}  t \in \mathcal{T}, j \in \mathcal{S}_t, i \in \mathcal{I}
\end{align}
\end{subequations}
We highlight that the proposed formulation can be extended in a number of ways to add different investment logic (e.g., account for limited investment budgets). Here, we present a formulation that contains enough features to highlight benefits of modular technologies in mitigating risk.

\section{Case Studies}
In this section, we present different case studies to illustrate how modularization can help mitigate risk. The first case study involves a single-product setting with 3 stages and has a structure of a binary tree. We then present a more complex and realistic case study that includes interdependent products and more stages and scenarios. The optimization problems were solved using {\em Gurobi} (version 0.7.6) with a default MIP Gap of $0.01\%$ and were implemented in the {\em JuMP} modeling framework. The scripts to reproduce all results can be found in \url{https://github.com/zavalab/JuliaBox/tree/master/ModularPlanning}.

\begin{figure}[hbt!]
\centering
    \includegraphics[scale=0.55]{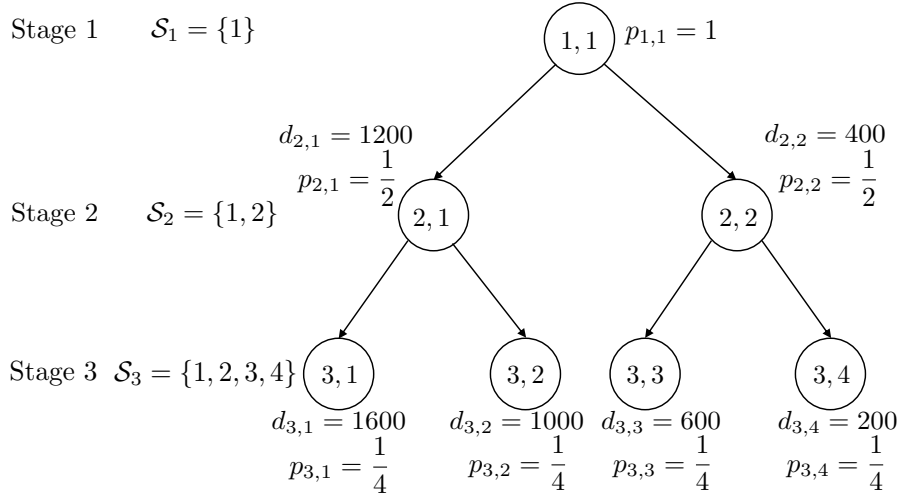}
    \caption{Tree representation for single-product stochastic case over 3 stages.}
    \label{fig:single}
\end{figure} 

\subsection{Single-Product Problem}

\begin{table}[!htb]
\caption {Data for single-product problem} \label{tab:single} 
\vspace{-0.2in}
\begin{center}
  \begin{tabular}{c*2{c}}
    \hline
Parameters & Values\\ \hline
Capacities (tons), $\mathcal{B}$ & \{100, 500, 1000, 1500\}\\ 
Installation Cost (\$), $\mathcal{C}$ & \{247, 721, 1145, 1500\}\\
Capacity Limit (tons), $\bar{x}$ & 1500\\
Installation Cost Limit (\$), $\bar{y}$ & 2000\\
Storage Cost (\$), $\rho_s$ & 30\\
Storage Limit (tons), $\bar{s}$ & 400\\
Waste Cost (\$), $\rho_w$ & 30\\
Operational Cost (\$), $\rho_p$ & 50\\
Selling Price (\$), $\pi_p$ & 140\\
Discount rate, $\gamma$ & 0.06\\\hline
  \end{tabular}
\end{center}
\end{table}

Figure \ref{fig:single} shows the stages, scenarios, and their corresponding demand and probabilities. The number inside each node represents the scenarios at each stage and the other number next to the node indicates the demand for each scenario in tons. Each parent node has two children nodes and we assume that each outcome has equal probability. All other required data is summarized in Table \ref{tab:single}. All the capacity-related quantities have the units of metric tons and price-related quantities have units of US dollars. 
\\

In this problem, we are seeking to make investment decisions at stage 1 and 2 that can help minimize NPV risk while achieving a constant level of expected NPV. To see the effect of that modular units have on flexibility, we solved this problem under three different capacity options (we call them Cases 1,2,3). In Case 1, we only allow the investor to choose between large capacities of 1500 tons and 1000 tons. In Case 2, we add medium capacity unit (500 tons ) to the list to provide more flexibility. In Case 3, we allow the investor to choose from the complete capacity list (which includes smaller modular units). The three cases are solved for the undiscounted and discounted NPV problem (to see the impact of time value of money). 
\vspace{0.1in}

\begin{figure}[!htp]
\centering
    \includegraphics[scale=0.52]{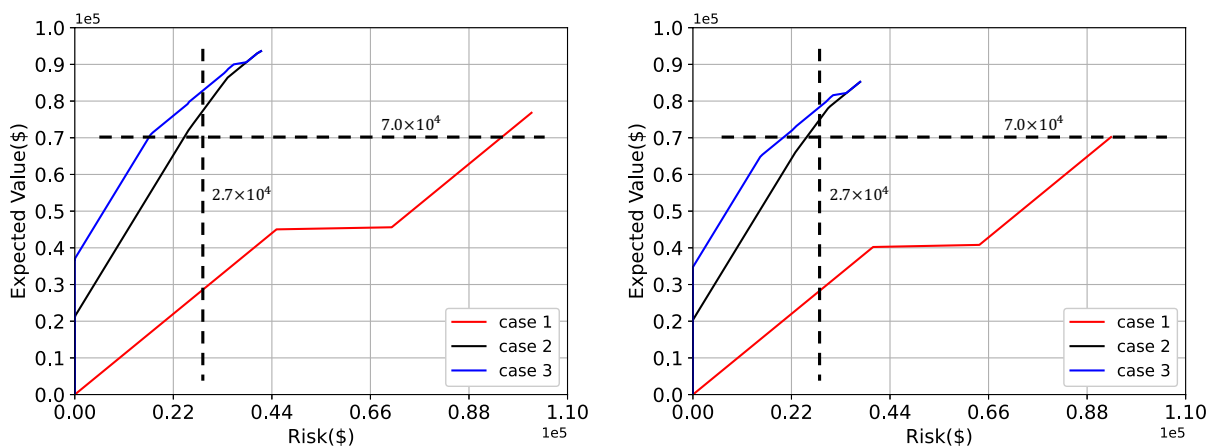}
    \caption{Pareto frontiers under undiscounted NPV (left) and discounted NPV (right) settings.}
    \label{fig:single_pareto}
\end{figure} 

The Pareto frontiers for both problems are shown in Figure \ref{fig:single_pareto}. Examples of investment plans obtained with these formulations are shown in Table \ref{tab:single_p}. We can see that the shape of the Pareto frontiers for the discounted and undiscounted problems is similar. We can thus see that the discounting factor does not influence the decisions made at each stage since the number of stages is small. As we will see in the next case study, the effect of discounting can be quite pronounced for problems that involve long planning horizons.  The Pareto frontiers highlight that that a strong trade-off exists between the expected value and risk of the NPV (higher expected NPV results in higher risk). This trade-off arises from economies of scale and flexibility (it is less expensive but more risky to install large units). It is clear that the Pareto frontier for cases 2 and 3 (under which small units are available) dominate the frontier of case 1 (under which only large units are available). Importantly, this occurs even if the installation costs of the large units have better economies of scale. At the same level for the expected NPV, Cases 2 and 3 achieve a significantly lower risk (reduction by a factor of 3). Similarly, at the same risk level, Cases 2 and 3 achieve a much higher expected NPV (increase by a factor of 2). We can also see that Case 2 and 3 achieve levels of expected NPV that are not achievable in Case 1.
\\

\begin{figure}[!htp]
\centering
    \includegraphics[scale=0.5]{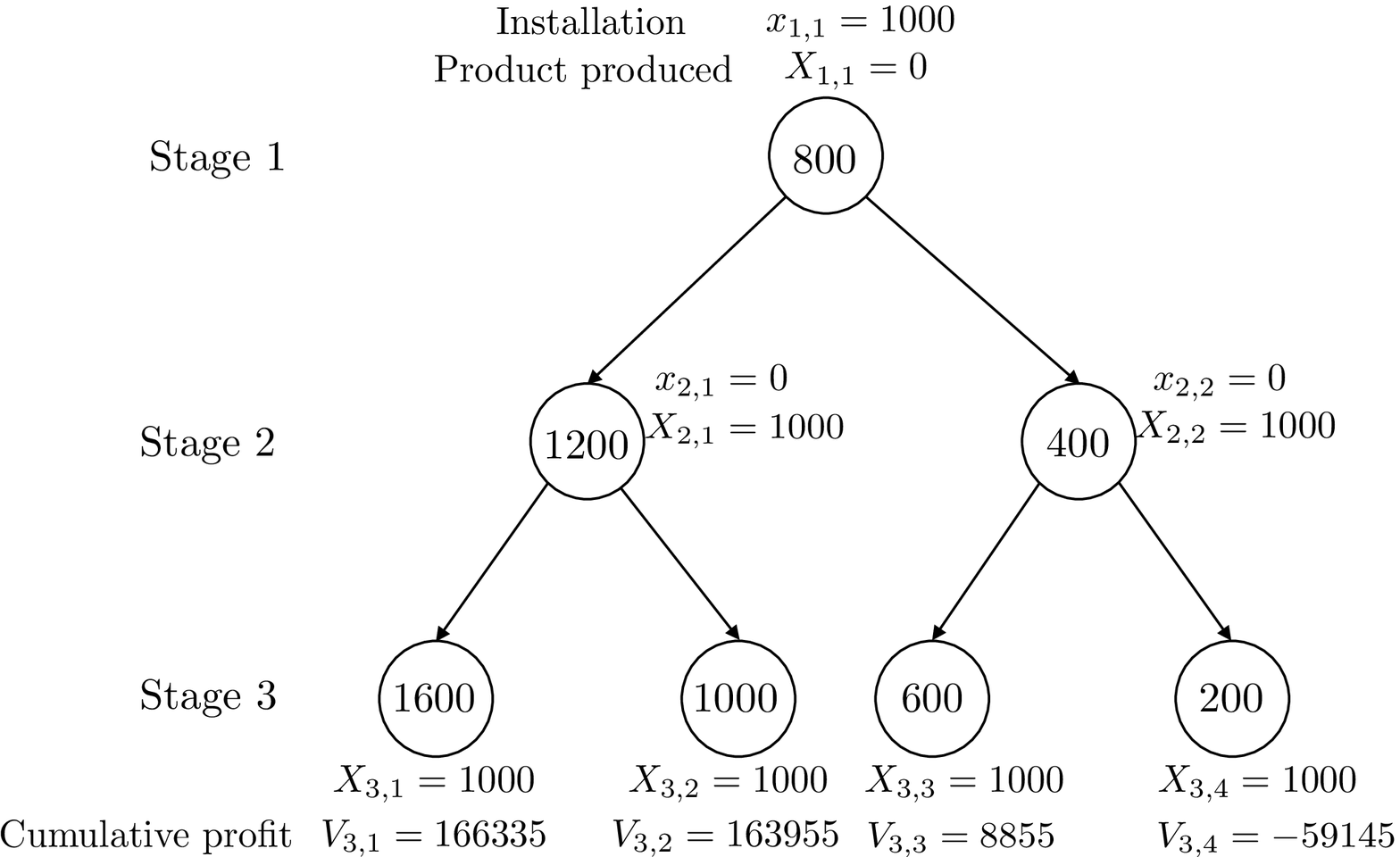}
    \caption{Investment strategy under Case 1 (undiscounted NPV setting).}
    \label{fig:single_result_1}
\end{figure}

\begin{figure}[!htp]
\centering
    \includegraphics[scale=0.5]{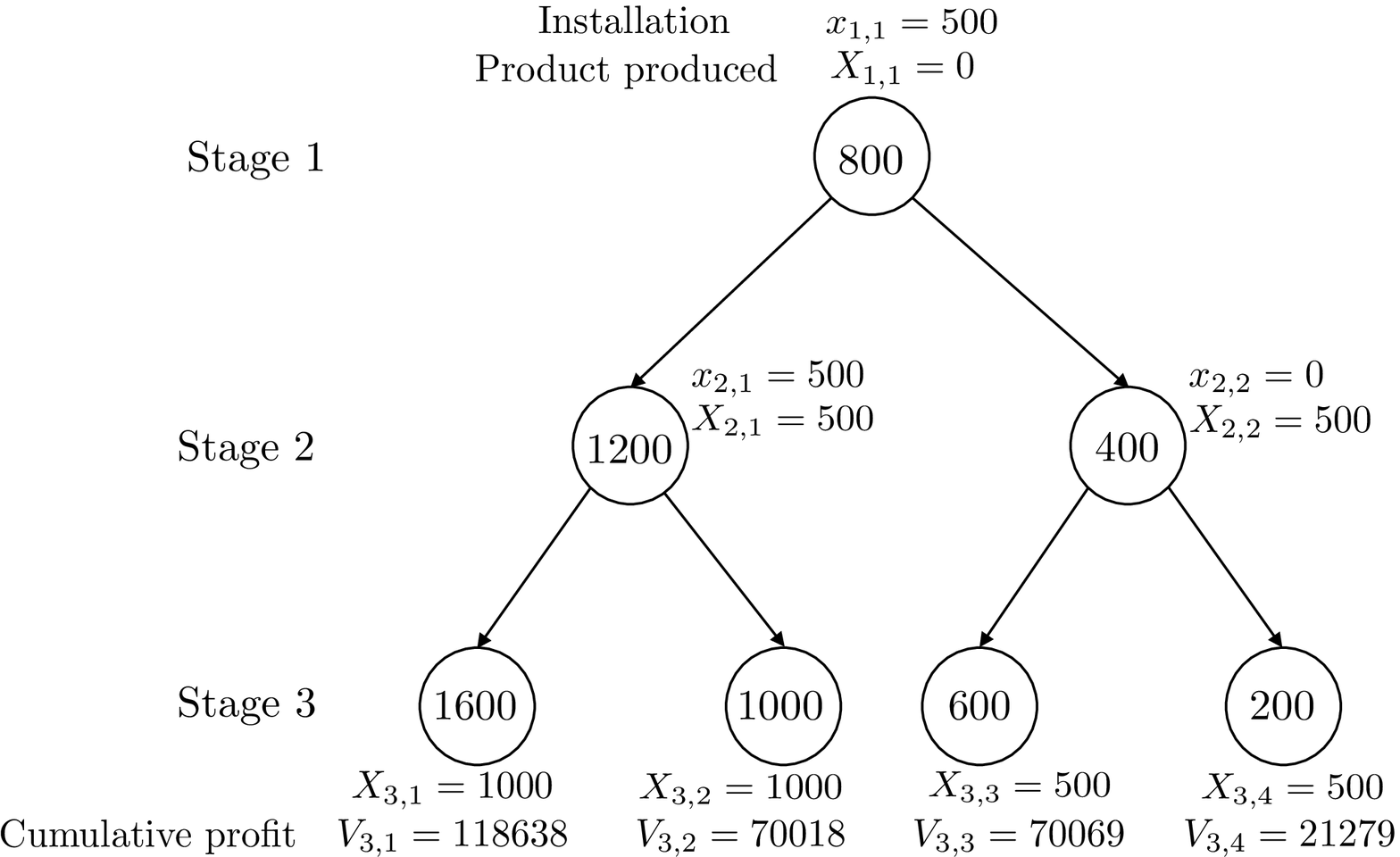}
    \caption{Investment strategy under Case 2 (undiscounted NPV setting).}
    \label{fig:single_result_2}
\end{figure}

\begin{figure}[!htp]
\centering
    \includegraphics[scale=0.5]{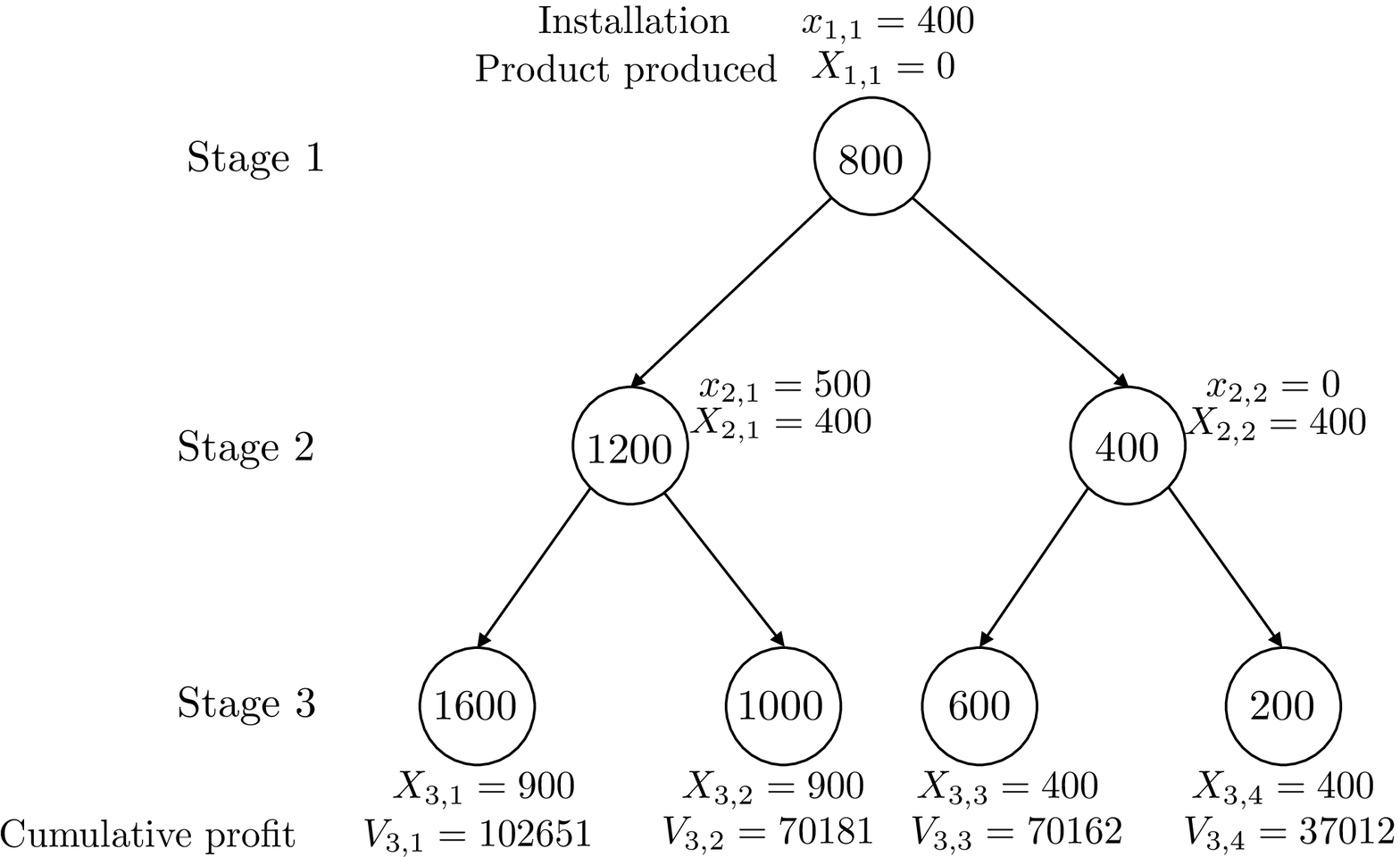}
    \caption{Investment strategy under Case 3 (undiscounted NPV setting).}
    \label{fig:single_result_3}
\end{figure}

In Table \ref{tab:single_p}, the installation column shows installation decisions made for the undiscounted problem. These decisions are also visualized in Figures \ref{fig:single_result_1}, \ref{fig:single_result_2}, and \ref{fig:single_result_3}. We can see that, with only larger capacity options (Case 1), we have no choice but to install the large unit at stage 1. In Case 2, we install a medium-sized unit in stage 1 and another medium-sized unit in stage 2; this achieves the same expected NPV but the risk is drastically reduced. For Case 3, we installed four small-sized units in the first stage and one medium-sized unit in the second stage. This achieves the same expected NPV but further decreases the risk. We thus conclude that the different capacity choices enable higher investment flexibility and reduced risk. 


\begin{table}[!htb]
\caption{Investment strategy and associated risks  (undiscounted NPV setting).}
\label{tab:single_p}
\vspace{-0.2in}
\begin{center}
  \begin{tabular}{c*3{c}}
    \hline
    Cases & Technology Sizes (tons) & Risk (\$) & Expected Value (\$)\\ \hline
    Case 1 & \{1000, 0, 0\} & 93149 & $7.0 \times 10^4$\\
    Case 2 & \{500, 500, 0\}&  24361 & $7.0 \times 10^4$\\
    Case 3 & \{100$\times$4, 500 , 0\} &  16495 & $7.0 \times 10^4$\\\hline
    \end{tabular}
\end{center}
\end{table}

\newpage

\begin{figure}[!htp]
\centering
    \includegraphics[scale=0.25]{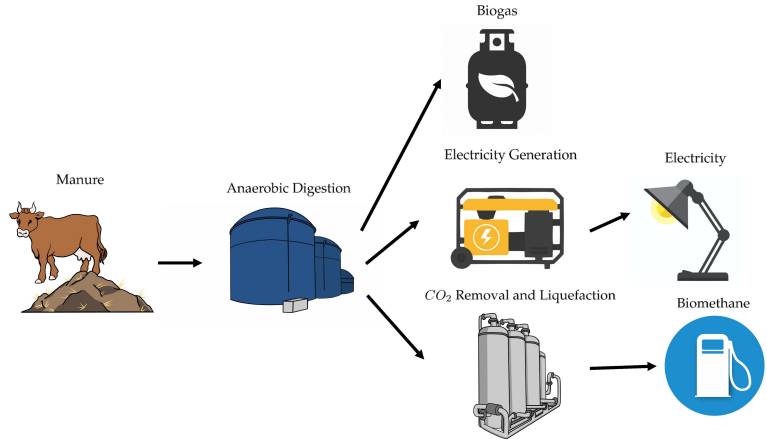}
    \caption{Process for the production of biogas and its byproducts}
    \label{fig:process_flow}
\end{figure} 

\begin{figure}[!htp]
\centering
    \includegraphics[scale=0.4]{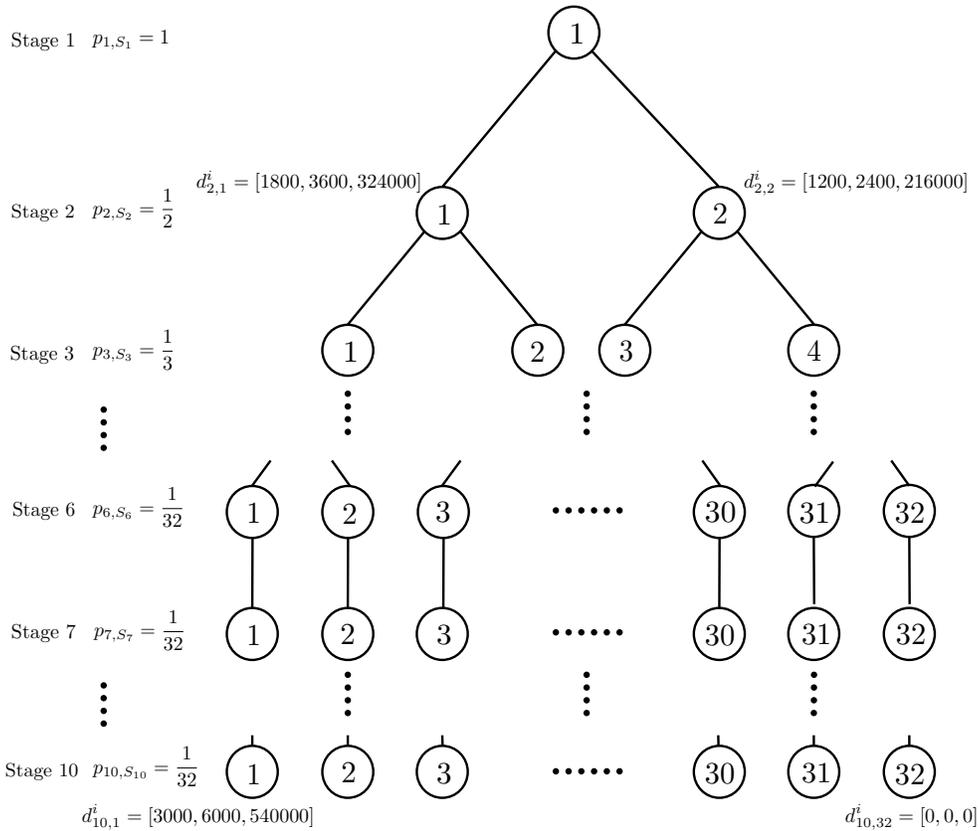}
    \caption{Tree representation of planning stages and scenarios of biogas case study.}
    \label{fig:multi}
\end{figure} 

\renewcommand{\arraystretch}{0.7}
\begin{table}[!htp]
\caption {Data for biogas capacity expansion problem.} \label{tab:multi} 
\vspace{-0.2in}
\begin{center}
  \begin{tabular}{c*2{c}}
    \hline
Parameters & Notation  & Values\\ \hline
\multirow{3}{*}{Capacities List }&$B^{i_1}$ (tons)&[400, 800, 1200]\\&$B^{i_2}$ (MWh)&[500, 1000, 2000]\\&$B^{i_3}$ (gallons)&[60000, 180000, 300000]\\ \hline
\multirow{3}{*}{Installation Cost (\$)}&$C^{i_1}$&[272844, 374693, 457138]\\&$C^{i_2}$&[172219, 234688, 347021]\\&$C^{i_3}$&[577269, 1930312, 2935611]\\\hline
\multirow{3}{*}{Capacity Limit}&$\bar{x}^{i_1}$ (tons)&6000\\&$\bar{x}^{i_2}$ (MWh)&6000\\&$\bar{x}^{i_3}$ (gallons)&1500000\\\hline
\multirow{3}{*}{Storage Cost}&$\rho_s^{i_1}$ (\$ per tons)&0\\&$\rho_s^{i_2}$ (\$ per MWh)&150000\\&$\rho_s^{i_3}$ (\$ per gallon)&0.2\\\hline
\multirow{3}{*}{Storage Limit}&$\bar{s}^{i_1}$ (tons)&800\\&$\bar{s}^{i_2}$ (MWh)&1000\\&$\bar{s}^{i_3}$ (gallons)&180000\\\hline
\multirow{3}{*}{Waste Disposal Cost}&$\rho_w^{i_1}$ (\$ per ton)&0\\&$\rho_w^{i_2}$ (\$ per MWh)&0\\&$\rho_w^{i_3}$ (\$ per gallon)&2\\\hline
\multirow{3}{*}{Operational Cost}&$\rho_p^{i_1}$ (\$ per ton)&34\\&$\rho_p^{i_2}$ (\$ per MWh)&40\\&$\rho_p^{i_3}$ (\$ per gallon)&0.56\\\hline
\multirow{3}{*}{Selling Price}&$\pi_p^{i_1}$ (\$ per ton)&100\\&$\pi_p^{i_2}$ (\$ per MWh)&130\\&$\pi_p^{i_3}$ (\$ per gallon)&2.5\\\hline
\multirow{2}{*}{Interdependency}&$\alpha_{i_2, i_1}$&0.68\\&$\alpha_{i_3, i_1}$&0.0046\\ \hline
Total Installation Cost Limit (\$) & $\bar{y}$ & $1 \times 10^7$\\ \hline
Interest Rate & $\gamma$ & 0.06 \\ \hline
  \end{tabular}
\end{center}
\end{table}

\subsection{Multi-Product Problem}

Biogas is a methane-rich gas mixture that can be produced from anaerobic digestion of organic waste (such as cow manure). The biogas (in metric tons) can be sold directly or can be used as raw material to produce electricity (in MWh) and liquefied biomethane (in gallons). These products are represented as $i_1$, $i_2$, and $i_3$ respectively \cite{Kirch_2005, Sampat_2018}. The process under study is visualized in Figure \ref{fig:process_flow}. The available technology capacities, investment cost, operation cost and other required information are summarized in Table \ref{tab:multi} \cite{Hu_2018, Beddoes_2007, energy_cost}. The capital cost for these technologies roughly follows the 2/3 rule.  The products are interdependent: producing 1 MWh of $i_2$ requires 0.68 tons of $i_1$ and producing 1 gallon of $i_3$ requires 0.0046 tons of $i_1$.  The planning stages have a duration of one year; as such, all capacities and production levels are expressed on a per-year basis. 
\\

The stochastic multistage setting is illustrated in Figure \ref{fig:multi}. Here, we have a planning horizon with 10 stages. From stage 1 to stage 6, each parent node has two children nodes (which capture variability in market demands); after stage 6, the market is assumed to stay constant and thus each parent node only has one children node. The demand for selected nodes is shown next to the node. For the first six stages, the children nodes of each parent node represent an optimistic market and a pessimistic market. 
\\

We consider 3 possible cases;  for Case 1, the unit for producing $i_1$ has a capacity of 1200 tons, the unit producing $i_2$ has a capacity of 2000 MWh, and the unit for producing $i_3$ has the capacity of 300000 gallons. In Case 2, we add a unit with a capacity of 800 tons for producing $i_1$, we add a unit with a capacity of 1000 MWh for the choices for $i_2$, and a capacity of 180,000 gallons for $i_3$. For Case 3, we further expand the capacity choices for product $i_1$ to include 400 tons, expand choices for $i_2$ to include 500 MWh, and expand choices for $i_3$ to include 60,000 gallons.  To provide some context on the size of these units, an annual capacity of 500 MWh corresponds to a power capacity of 500/8760=0.057 MW (57 kW). As such, the small capacities for the power generators correspond to those of small modular systems. 

\begin{figure}[!htp]
\centering
    \includegraphics[scale=0.52]{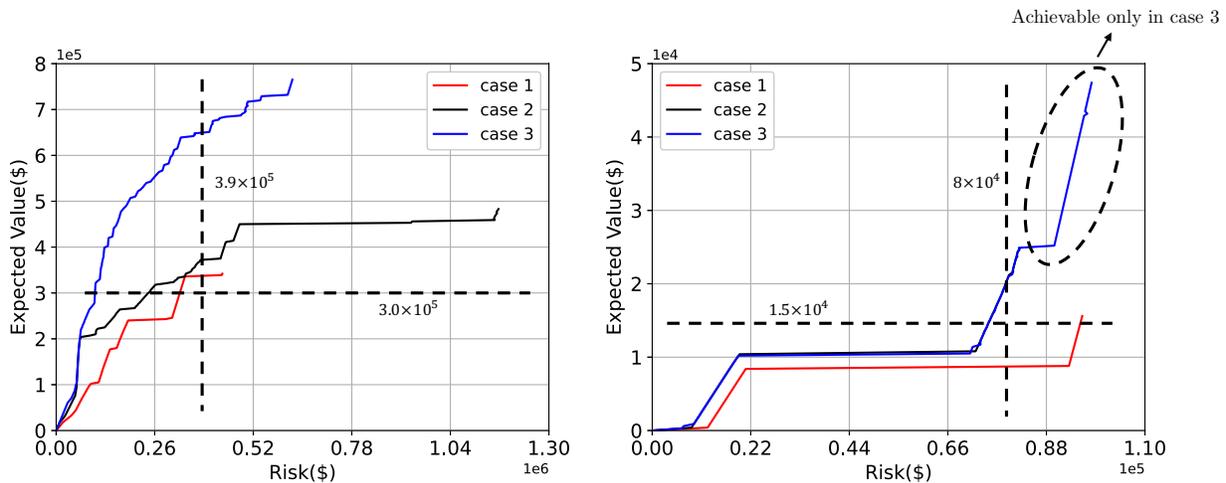}
    \caption{Pareto frontier for undiscounted NPV (left) and discounted NPV (right) settings.}
    \label{fig:multi_pareto}
\end{figure} 

We would like to determine the number of planning stages that it takes for the investment to be profitable. As such, we gradually increase the planning horizon of the CE problem until the profit is positive. We found that, for both discounted and undiscounted problems, the expected NPV remains zero for any planning horizon with less than 8 stages. In other words, the investment is only profitable if the project lifetime is at least 8 years. We can thus see that the length of the planning horizon plays an important role in making investment decisions. We assume that the planning horizon is 10 years (as shown in Figure \ref{fig:multi}).  Again, we would like to determine an optimal investment strategy that maximizes expected NPV and minimizes its risk. The Pareto frontiers are shown in Figure \ref{fig:multi_pareto} and we compare risks obtained under the different cases in Table \ref{tab:multi_1} and Table \ref{tab:multi_2}.

\setlength{\tabcolsep}{3pt}
\renewcommand{\arraystretch}{0.8}
\begin{table}[!htb]
\caption {Investment strategy for undiscounted NPV problem with $\mathcal{E}=3.0 \times 10^5$} \label{tab:multi_1} 
\vspace{-0.2in}
\begin{center}
  \begin{tabular}{c*6{c}}
    \hline
\multirow{2}{*}{Cases} & \multirow{2}{*}{Installation} & \multirow{2}{*}{Expected Value (\$)} & \multirow{2}{*}{Risk (\$)} & \multirow{2}{*}{ \# of Constraints} & \# of Variables\\&&&&& (Cont.+Int.)\\ \hline
\multirow{2}{*}{Case 1}&$x_{1,1}^{i_1}=1200 \times 2$ tons&\multirow{2}{*}{$3.0 \times 10^5$}&\multirow{2}{*}{$3.26 \times 10^5$ }& \multirow{2}{*}{4355}  & \multirow{2}{*}{2225+1623} \\&$x_{1,1}^{i_2}=2000$ MWh\\ \hline
\multirow{8}{*}{Case 2}&$x_{1,1}^{i_1}=1200$ tons&\multirow{8}{*}{$3.0 \times 10^5$}&\multirow{8}{*}{$2.43 \times 10^5$ }& \multirow{8}{*}{4355}  & \multirow{8}{*}{2225+2100} \\&$x_{1,1}^{i_2}=1000$ MWh\\&$x_{2,1}^{i_1}=1200$ tons\\&$x_{2,1}^{i_3}=180000$ gallons\\&$x_{3,1}^{i_1}=1200$ tons\\&$x_{5,1}^{i_2}=2000$ MWh\\&$x_{5,4}^{i_2}=1000$ MWh\\&$x_{6,6}^{i_2}=2000$ MWh\\\hline
\multirow{4}{*}{Case 3}&$x_{1,1}^{i_1}=800$ tons&\multirow{4}{*}{$3.0 \times 10^5$}&\multirow{4}{*}{$1.02 \times 10^5$ }& \multirow{4}{*}{4355}  & \multirow{4}{*}{2225+2577}\\&$x_{1,1}^{i_3}=60000$ gallons\\&$x_{3,1}^{i_1}=1200 \times 2$ tons\\&$x_{3,1}^{i_3}=60000 \times 2$ gallons\\\hline
  \end{tabular}
\end{center}
\end{table}

\setlength{\tabcolsep}{8pt}
\renewcommand{\arraystretch}{0.8}
\begin{table}[!htb]
\caption {Investment strategy for discounted NPV problem with $\mathcal{R}=3.9 \times 10^5$} \label{tab:multi_3} 
\vspace{-0.2in}
\begin{center}
  \begin{tabular}{c*4{c}}
    \hline
Cases & Installation & Expected Value(\$) & Risk (\$)\\ \hline
\multirow{2}{*}{Case 1}&$x_{1,1}^{i_1}=1200 \times 2$ tons&\multirow{2}{*}{$2.64 \times 10^5$}&\multirow{2}{*}{$3.90 \times 10^5$ }\\&$x_{1,1}^{i_2}=2000$ MWh\\ \hline
\multirow{6}{*}{Case 2}&$x_{1,1}^{i_1}=1200$ tons&\multirow{6}{*}{$3.72 \times 10^5$}&\multirow{6}{*}{$3.90 \times 10^5$ }\\&$x_{2,1}^{i_1}=1200 \times 2$ tons\\&$x_{2,1}^{i_2}=2000$ MWh\\&$x_{2,1}^{i_3}=180000$ gallons\\&$x_{3,4}^{i_2}=1000$ MWh\\&$x_{6,24}^{i_2}=1000$ MWh\\\hline
\multirow{5}{*}{Case 3}&$x_{1,1}^{i_1}=1200$ tons&\multirow{5}{*}{$6.50 \times 10^5$}&\multirow{5}{*}{$3.90 \times 10^5$ } \\&$x_{1,1}^{i_3}=60000 \times 2$ gallons\\&$x_{2,1}^{i_1}=1200 \times 2$ tons\\&$x_{2,1}^{i_2}=2000$ MWh \\ &$x_{3,1}^{i_3}=60000$ tons\\\hline
  \end{tabular}
\end{center}
\end{table}

\setlength{\tabcolsep}{8pt}
\renewcommand{\arraystretch}{0.8}
\begin{table}[!htb]
\caption {Investment strategy for discounted NPV problem with $\mathcal{E}=1.5 \times 10^4$} \label{tab:multi_2} 
\vspace{-0.2in}
\begin{center}
  \begin{tabular}{c*4{c}}
    \hline
Cases & Installation & Expected Value(\$) & Risk (\$)\\ \hline
Case 1&$x_{1,1}^{i_1}=1200$ tons&$1.5 \times 10^4$&$3.06 \times 10^6$\\ \hline
\multirow{4}{*}{Case 2}&$x_{1,1}^{i_1}=1200$ tons&\multirow{4}{*}{$1.5 \times 10^4$}&\multirow{4}{*}{$2.41 \times 10^6$ }\\&$x_{3,4}^{i_2}=1000$ MWh\\&$x_{5,12}^{i_2}=1000$ MWh\\&$x_{6,16}^{i_2}=1000$ MWh\\\hline
\multirow{4}{*}{Case 3}&$x_{1,1}^{i_1}=1200$ tons&\multirow{4}{*}{$1.5 \times 10^4$}&\multirow{4}{*}{$2.41 \times 10^6$ } \\&$x_{3,4}^{i_2}=1000$ MWh\\&$x_{5,12}^{i_2}=1000$ MWh\\&$x_{6,16}^{i_2}=1000$ MWh \\\hline
  \end{tabular}
\end{center}
\end{table}

\vspace{0.1in}
We again find that the Pareto frontier of Case 3 (considering small technologies) dominates. Looking horizontally  (for the same expected NPV) cases with more capacity options reduce risk. Looking vertically (for the same risk) we can see that modularity allows us to reach higher expected NPVs. For the discounted NPV problem we see that adding smaller capacity options (Case 2) reduces risk, further reducing the capacity (Case 3) can achieve higher expected profits but does not help to mitigate the risk. This is because of complex interplays between discounting and economies of scale. As we discount the future cash flow, the effect of installing small capacities at future stages reduces, and together with the effect of economies of scale, the advantages brought by modular technologies become less obvious. This indicates that reducing technology sizes aids flexibility (but there is a limit to such flexibility). From Table \ref{tab:multi_1} and Table \ref{tab:multi_3}  we can see that, for the undiscounted problem, most of the investment occurs at the early stages. Here, we can also see that modular technologies are used extensively to reduce risk (risk is reduced by a factor of three) and increase profit (profit is increased by a factor of three).

\section{Conclusion and Future Work}
We examined the ability of using small modular technologies to control investment risk. To do so, we propose a capacity expansion problem that aims to determine optimal investment strategies over a given planning horizon. This expansion problem is a stochastic, multistage, and multiobjective optimization problem. We propose to measure joint risk on the net present value over the entire planning horizon in order to avoid the ambiguity associated with standard multistage stochastic formulations. Our analysis reveals that small technologies provide flexibility that translates into tangible reductions of risk (despite the fact that they are not benefited by economies of scale). However, we also find that flexibility provided by capacity reductions has limits that result from the complex interplay between economies of scale and discounting. As part of future work, we are interested in exploring the use of decomposition strategies to address tractability issues (e.g., by using stochastic dual dynamic programming techniques). In this work we ignored engineering costs associated with different types of technologies (which can be reduced using modularization). We will use more detailed cost representations and case studies in future work.

\section*{Acknowledgments}
We acknowledge funding from the NSF CAREER award CBET-1748516.
We also acknowledge partial support from the U.S. Department of Agriculture  under grant 2017-67003-26055.
\bibliographystyle{abbrv}
\bibliography{ref}

\end{document}